\documentclass[11pt]{article}
\usepackage{amsthm}
\usepackage{pxfonts}
\usepackage{yfonts}
\usepackage{dsfont}
\usepackage{graphicx}
\usepackage{relsize}
\usepackage{color}

\def\bel{\begin{equation}\label}
\def\eeq{\end{equation}}
\def\ds{\displaystyle}

\def\mt{\longrightarrow}
\def\v{\vskip 1em}

\def\R{\mathbb R}

\def\C{\mathfrak{B}}

\def\N{{\bf N}}

\def\Re{{\bf Re}}

\def\S{{\bf S}}

\def\F{{\bf F}}

\def\L{{\bf L}}

\def\V{{\bf V}}

\def\T{{\bf T}}
\def\X{{\bf X}}
\def\Y{{\bf Y}}

\def\p{{\partial}}
\def\a{{\bf a}}
\def\b{{\bf b}}

\def\bar{\overline}

\def\I{{\bf I}}

\def\alpha{\alphaup}
\def\beta{\betaup}
\def\gamma{\gammaup}
\def\delta{\deltaup}
\def\theta{\thetaup}
\def\xi{{\xiup}}
\def\eta{{\etaup}}
\def\tau{{\tauup}}
\def\rho{{\rhoup}}
\def\phi{{\phiup}}
\def\psi{{\psiup}}
\def\lambda{{\lambdaup}}
\def\omega{\omegaup}
\def\varphi{{\varphiup}}
\def\gamma{{\gammaup}}
\def\c{{\bf c}}

\newtheorem{remark}{Remark}[section]

\begin{document}

\[\hbox{\LARGE{\bf Hardy-Littlewood-Sobolev inequality revisit}}\]
\[\hbox{\LARGE{\bf on Heisenberg group}}\]

\[\hbox{Chuhan Sun~~~~and~~~~Zipeng Wang}\]
\begin{abstract}
We study a family of fractional integral operators  
\[\I_{\alpha\beta}f(u,v,t)~=~\iiint_{\R^{2n+1}} f(\xi,\eta,\tau)\V^{\alpha\beta}\big[(u,v,t)\odot(\xi,\eta,\tau)^{-1}\big]d\xi d\eta d\tau\]
where $(u,v,t)\odot(\xi,\eta,\tau)^{-1}=\big[u-\xi, v-\eta,t-\tau-\mu(u\cdot\eta-v\cdot\xi)\big]$, $\mu\in\R$.
\v
$\V^{\alpha\beta}$ is a distribution in $\R^{2n+1}$ satisfying Zygmund dilations.
A characterization is established between the $\L^p\mt\L^q$-boundedness of $\I_{\alpha\beta}$ and the necessary constraints consisting of $\alpha,\beta\in\R$ and $1<p<q<\infty$. 
\end{abstract}

\section{Introduction}
\setcounter{equation}{0}
Let $0<\a<\N$. A fractional integral operator $\T_\a$ is initially defined on $\R^\N$ as
\bel{I_a}
\T_\a f(x)~=~\int_{\R^\N}f(y) \left[{1\over|x-y|}\right]^{\N-\a} dy.
\eeq
In 1928, Hardy and Littlewood \cite{Hardy-Littlewood} have obtained an regularity theorem for $\T_\a$ when $\N=1$. Ten years later, Sobolev \cite{Sobolev} made extensions on every higher dimensional space.

\v
{\bf Hardy-Littlewood-Sobolev theorem}~~~{\it Let $\T_\a$ defined in (\ref{I_a}) for $0<\a<\N$. We have
\bel{HLS Ineq}
\begin{array}{cc}\ds
\left\| \T_\a f\right\|_{\L^q(\R^\N)}~\leq~\C_{p~q}~\left\| f\right\|_{\L^p(\R^\N)},\qquad 1<p<q<\infty
\\\\ \ds
\hbox{\small{if and only if}}\qquad {\a\over \N}~=~{1\over p}-{1\over q}.~~~~~~~
\end{array}
\eeq
}
$\diamond$ {\small Throughout, $\C>0$ is regarded as a generic constant depending on its sub-indices.}

This classical result was first re-investigated by Folland and Stein \cite{Folland-Stein} on Heisenberg group.
We shall be working on its real variable representation with a multiplication law:
\bel{multiplication law}
\begin{array}{cc}\ds
(u,v,t)\odot(\xi,\eta,\tau)~=~\Big[u+\xi, v+\eta,t+\tau+\mu\big(u\cdot\eta-v\cdot\xi\big)\Big],
\\\\ \ds
(u,v,t)~\in~\R^n\times\R^n\times\R,\qquad (\xi,\eta,\tau)^{-1}~=~(-\xi,-\eta,-\tau)~\in~\R^n\times\R^n\times\R.
\end{array}
\eeq
whenever $\mu\in\R$.

Let $0<\a<n+1$. Consider
\bel{T_a}
\S_\a f(u,v,t)~=~\iiint_{\R^{2n+1}}f(\xi,\eta,\tau)\Omega^\a\Big[(u,v,t)\odot(\xi,\eta,\tau)^{-1}\Big]d\xi d\eta d\tau.
\eeq
$\Omega^\a$ is a distribution in $\R^{2n+1}$ agree with
\bel{Omega^a}
\Omega^\a(\xi,\eta,\tau)~=~\left[{1\over |\xi|^2+|\eta|^2+|\tau|}\right]^{n+1-\a},\qquad \hbox{\small{$(\xi,\eta,\tau)\neq(0,0,0)$}}.
\eeq
Observe that
\bel{non-isotropic}
\Omega^\a\Big[(\delta u,\delta v,\delta^2 t)\odot(\delta \xi,\delta \eta,\delta^2\tau)^{-1}\Big]~=~\delta^{2\a-2n-2}\Omega^\a\Big[(u,v,t)\odot(\xi,\eta,\tau)^{-1}\Big], \qquad \delta>0.
\eeq

{\bf Folland-Stein theorem} ~~{\it Let $\S_\a$ defined in (\ref{T_a})-(\ref{Omega^a}) for $0<\a<n+1$. We have
\bel{Folland-Stein theorem}
\begin{array}{cc}\ds
\left\| \S_\a f\right\|_{\L^q(\R^{n+1})}~\leq~\C_{p~q}~\left\| f\right\|_{\L^p(\R^{2n+1})},\qquad 1<p<q<\infty
\\\\ \ds
\hbox{if and only if}\qquad {\a\over n+1}~=~{1\over p}-{1\over q}.
\end{array}
\eeq}\\
The best constant for the $\L^p\mt\L^q$-norm inequality in (\ref{Folland-Stein theorem}) is found by Frank and Lieb \cite{Frank-Lieb}. A discrete analogue of this result has been obtained by Pierce \cite{Pierce}. Recently, the regarding commutator estimates are established  by Fanelli and~Roncal \cite{Luca-Luz}.

In this paper, we introduce a family of fractional integral operators  whose kernels have a mixture of homogeneities defined in $\R^{2n+1}$ with a multiplication law $\odot$ in (\ref{multiplication law}).  An initial motivation for considering such operators that commute with multi-parameter dilations  comes from the $\bar{\p}$-Neumann problem on the model domain which has a Heisenberg group as its boundary.  The unique solution turns out to be a composition of two singular integral operators. One of them is elliptic associated with a standard one-parameter dilation. The other is parabolic whose kernel satisfies an non-isotropic dilation as (\ref{non-isotropic}). Singular integrals of this type have been systematically studied by Phong and Stein \cite{Phong-Stein} and later refined by 
Muller, Ricci and Stein \cite{Muller-Ricci-Stein}. 

One particularly interesting example among certain operators having a negative order is  $\mathcal{L}^{-\a}T^{-\b}$ for $0<\a<n, 0<\b<1$ and $\a\ge n\b$ where $T=\partial_t$ and $\mathcal{L}$ is the sub-Laplacian:
$
\mathcal{L}=-\sum_{j=1}^n \X_j^2+\Y^2_j,~ \X_j=\p_{x_j}+2 y_j \p_t,~ \Y_j=\p_{y_j}-2x_j \p_t$.
The inverse of $\mathcal{L}^\a$, $\Re\a>0$  is given as the Riesz potential defined on Heisenberg group. 
Namely,  
\[
\mathcal{L}^{-\a}~=~{1\over \Gamma(\a)}\int_0^\infty s^{\a-1} e^{-s\mathcal{L}} ds,\qquad\hbox{\small{$\Re\a>0$}}
\]
where $\Gamma$ is Gamma function. More background can be found in chapter XIII of Stein \cite{Stein}.

Let $0<\a<n$, $0<\b<1$ and $\a\ge n\b$. We have
\bel{Muller-Ricci-Stein Result}
\begin{array}{cc}\ds
\left\|\mathcal{L}^{-\a}T^{-\b} f\right\|_{\L^q(\R^{2n+1})}~\leq~\C_{\a~\b~p}~\left\| f\right\|_{\L^p(\R^{2n+1})},\qquad 1<p<q<\infty
\\\\ \ds
\hbox{\small{if and only if}}\qquad {\a+\b\over n+1}~=~{1\over p}-{1\over q}.
\end{array}
\eeq
This $\L^p\mt\L^q$-regularity result  is proved by using complex interpolation in section 6 of \cite{Muller-Ricci-Stein}. One of the two end-point estimates relies on the $\L^p$-theorem developed thereby.

Let $0<\a<n$, $0<\b<1$ and $\a\ge n\b$. $\Omega^{\a\b}$ is a distribution in $\R^{2n+1}$ agree with 
\bel{Omega^ab}
\Omega^{\a\b}(\xi,\eta,\tau)~=~\left[{1\over |\xi|^2+|\eta|^2}\right]^{n-\a}\left[{1\over |\xi|^2+|\eta|^2+|\tau|}\right]^{1-\b},\qquad\hbox{\small{$(\xi,\eta)\neq(0,0)$}}.
\eeq
The kernel of $\mathcal{L}^{-\a}T^{-\b}$ is similar to $\Gamma\Big({1-\b\over 2}\Big)\Omega^{\a\b}(\xi,\eta,\tau)$ for  $(\xi,\eta)\neq(0,0)$. See {\bf Theorem 6.2} in \cite{Muller-Ricci-Stein}.  (We say $A$ similar to $B$ if $\c^{-1} B\leq A\leq \c B$ for some $\c>0$.)

{\bf Question}: 
{\it For every operator having a kernel  similar to $\Omega^{\a\b}$ away from its singularity, does it satisfy the regularity estimate in (\ref{Muller-Ricci-Stein Result})?}

 The answer  is yes. Consider
\bel{S_ab}
\S_{\a\b}f(u,v,t)~=~\iiint_{\R^{2n+1}} f\left(\xi,\eta,\tau\right) \Omega^{\a\b}\Big[(u,v,t)\odot(\xi,\eta,\tau)^{-1}\Big]d\xi d\eta d\tau.
\eeq

{\bf Theorem One}~~{\it Let  $\S_{\a\b}$ defined in (\ref{Omega^ab})-(\ref{S_ab}) for $0<\a<n$, $0<\b<1$ and $\a\ge n\b$. We have
\bel{Result One}
\begin{array}{cc}\ds
\left\| \S_{\a\b} f\right\|_{\L^q(\R^{2n+1})}~\leq~\C_{\a~\b~p}~\left\| f\right\|_{\L^p(\R^{2n+1})},\qquad 1<p<q<\infty
\\\\ \ds
\hbox{\small{if and only if}}\qquad {\a+\b\over n+1}~=~{1\over p}-{1\over q}.
\end{array}
\eeq}
\begin{remark} $\a\ge n \b$ is in fact an necessary condition for (\ref{Result One}).
\end{remark}
Let $\alpha, \beta\in\R$. $\V^{\alpha\beta}$ is a distribution in $\R^{2n+1}$ agree with 
\bel{V}
\V^{\alpha\beta}(\xi,\eta,\tau)~=~|\xi|^{\alpha-n}|\eta|^{\alpha-n}|\tau|^{\beta-1} \Bigg[ {|\xi||\eta|\over |\tau|}+{|\tau|\over |\xi||\eta|}\Bigg]^{-{|\alpha-n\beta |\over n+1}},\qquad\hbox{\small{$\xi\neq0, \eta\neq0$, $\tau\neq0$}}.
\eeq
\begin{remark} $\Omega^{\a\b}(\xi,\eta,\tau)$  in (\ref{Omega^ab}) can be bounded by two $\V^{\alpha\beta}(\xi,\eta,\tau)$ in (\ref{V}) for some $\alpha,\beta\in\R$ and $\a+\b=\alpha+\beta$.  
\end{remark}
Define 
\bel{I}
\begin{array}{lr}\ds
\I_{\alpha\beta} f(u,v,t)~=~\iiint_{\R^{2n+1}} f\left(\xi,\eta,\tau\right) \V^{\alpha\beta}\Big[(u,v,t)\odot(\xi,\eta,\tau)^{-1}\Big]d\xi d\eta d\tau.
\end{array}
\eeq
Observe that 
\bel{Zygmund dilations}
\begin{array}{lr}\ds
\V^{\alpha\beta}\Big[(\delta_1u,\delta_2v,\delta_1\delta_2t)\odot(\delta_1\xi,\delta_2\eta,\delta_1\delta_2\tau)^{-1}\Big]~=~\delta_1^{\alpha+\beta-n-1}\delta_2^{\alpha+\beta-n-1} \V^{\alpha\beta}\Big[(u,v,t)\odot(\xi,\eta,\tau)^{-1}\Big],
\\ \ds \delta_1,\delta_2>0.
\end{array}
\eeq
The two-parameter dilation in (\ref{Zygmund dilations}) is an example of Zygmund dilations:   $(u,v,t)\mt(\delta_1 u,\delta_2 v, \delta_1\delta_2 t), \delta_1,\delta_2>0$.
About maximal functions and singular integrals associated with Zygmund dilations, a number of pioneering results have been accomplished. For instance, see  Nagel and Wainger \cite{Nagel-Wainger}, Ricci and Stein \cite{Ricci-Stein}, Fefferman and Pipher \cite{Fefferman-Pipher}, Han et-al \cite{Han-} and Hytonen et-al \cite{Hytonen-}.
The area remains largely open for fractional integration. Our main result is stated in below.

{\bf Theorem Two}~~{\it Let $\I_{\alpha\beta}$ defined in  (\ref{V})-(\ref{I}) for $\alpha,\beta\in\R$. We have
\bel{Result Two}
\begin{array}{cc}\ds
\left\| \I_{\alpha\beta} f\right\|_{\L^q(\R^{2n+1})}~\leq~\C_{p~q}~\left\| f\right\|_{\L^p(\R^{2n+1})},\qquad 1<p<q<\infty
\\\\ \ds
\hbox{if and only if}\qquad
{\alpha+\beta\over n+1}~=~{1\over p}-{1\over q}.
\end{array}
\eeq}
\begin{remark} ${|\alpha- n\beta|\over n+1}$ given in (\ref{V}) is the smallest (best) exponent for which we can have  (\ref{Result Two}).
\end{remark} 
{\bf Theorem Two} implies {\bf Theorem One} because of {\bf Remark 1.2}. The rest of paper is organized as follows. In section 2, we prove some necessary constraints consisting of  $\a,\b$, $\alpha,\beta$ and $p,q$.  
These include {\bf Remark 1.1}, {\bf Remark 1.3} and the homogeneity condition in (\ref{Result One}) and (\ref{Result Two}). In section 3, we show {\bf Remark 1.2}. In section 4, we prove {\bf Theorem Two}.

\section{Some necessary constraints}
\setcounter{equation}{0}
Let $\S_{\a\b}$ defined in (\ref{Omega^ab})-(\ref{S_ab}) for $0<\a<n, 0<\b<1$ and $f\ge0$. By changing variable $\tau\mt \tau+\mu\big(u\cdot\eta-v\cdot\xi\big)$, we find
\bel{S rewrite}
\begin{array}{lr}\ds
\S_{\a\b} f(u,v,t)~=~
\\ \ds
\iiint_{\R^{2n+1}} f\left(\xi,\eta,\tau+\mu\big(u\cdot\eta-v\cdot\xi\big)\right) 
\left[{1\over |u-\xi|^2+|v-\eta|^2}\right]^{n-\a}\left[{1\over |u-\xi|^2+|v-\eta|^2+|t-\tau|}\right]^{1-\b}
\\ \ds~~~~~~~~~~~~~~~~~~~~~~~~~~~~~~~~~~~~~~~~~~~~~~~~~~~~~~~~~~~~~~~~~~~~~~~~~~~~~~~~~~~~~~~~~~~~~~~~~~~~~~~~~~~~~~~~~~~~~~~~~~~~~~~~~~~~
d\xi d\eta d\tau.
\end{array}
\eeq
By changing dilations $(u,v,t)\mt (\delta u, \delta v, \delta^2\lambda t)$ and $(\xi,\eta,\tau)\mt (\delta \xi, \delta \eta, \delta^2\lambda \tau)$ for $\delta>0, \lambda>1$,  we have
\bel{Dila S}
\begin{array}{lr}\ds
\left\{ \iiint_{\R^{2n+1}} 
\left\{\iiint_{\R^{2n+1}} f\left[\delta^{-1} \xi,\delta^{-1}\eta,\delta^{-2}\lambda^{-1}\big[\tau+\mu\lambda\big(u\cdot\eta-v\cdot\xi\big)\big]\right] \right.\right.
\\ \ds
\left. \left.\left[{1\over |u-\xi|^2+|v-\eta|^2}\right]^{n-\a}\left[{1\over |u-\xi|^2+|v-\eta|^2+|t-\tau|}\right]^{1-\b}
 d\xi d\eta d\tau\right\}^q dudvdt\right\}^{1\over q}
\\\\ \ds
=~\delta^{2\big[\a+\b\big]}\delta^{2n+2\over q}\lambda^{1\over q}~\left\{ \iiint_{\R^{2n+1}} 
\left\{\iiint_{\R^{2n+1}} f\left( \xi,\eta,\tau+\mu\big(u\cdot\eta-v\cdot\xi\big)\right) \right.\right.
\\ \ds~~~
\left. \left.\left[{1\over |u-\xi|^2+|v-\eta|^2}\right]^{n-\a}\left[{1\over |u-\xi|^2+|v-\eta|^2+\lambda|t-\tau|}\right]^{1-\b}
\lambda d\xi d\eta d\tau\right\}^q dudvdt\right\}^{1\over q}
\\\\ \ds
\ge~\delta^{2\big[\a+\b\big]}\delta^{2n+2\over q}\lambda^\b\lambda^{1\over q}~\left\{ \iiint_{\R^{2n+1}} 
\left\{\iiint_{\R^{2n+1}} f\left(\xi,\eta,\tau+\mu\big(u\cdot\eta-v\cdot\xi\big)\right) \right.\right.
\\ \ds~~~
\left. \left.\left[{1\over |u-\xi|^2+|v-\eta|^2}\right]^{n-\a}\left[{1\over |u-\xi|^2+|v-\eta|^2+|t-\tau|}\right]^{1-\b}
 d\xi d\eta d\tau\right\}^q dudvdt\right\}^{1\over q}.
\end{array}
\eeq
The $\L^p\mt\L^q$-norm inequality in (\ref{Result One}) implies that the last line of (\ref{Dila S}) is bounded by
\bel{S L^p}
\begin{array}{lr}\ds
\left\{ \iiint_{\R^{2n+1}} \Big[ f\left(\delta^{-1} \xi,\delta^{-1}\eta,\delta^{-2}\lambda^{-1} \tau\right) \Big]^p d\xi d\eta d\tau\right\}^{1\over p}~=~\delta^{2n+2\over p} \lambda^{1\over p} \left\| f\right\|_{\L^p(\R^{2n+1})}.
\end{array}
\eeq
This must be true for every $\delta>0$ and $\lambda>1$. We necessarily have
\bel{Constraints ab}
{\a+\b\over n+1}~=~{1\over p}-{1\over q},\qquad \b~\leq~{1\over p}-{1\over q}
\eeq
which together imply $\a\ge n\b$.

Let $\I_{\alpha\beta}$ defined in (\ref{V})-(\ref{I}) for $\alpha,\beta\in\R$ and $f\ge0$. By changing variable $\tau\mt \tau+\mu\big(u\cdot\eta-v\cdot\xi\big)$, we find
\bel{I rewrite}
\begin{array}{lr}\ds
\I_{\alpha\beta} f(u,v,t)~=~\iiint_{\R^{2n+1}} f\left(\xi,\eta,\tau+\mu\big(u\cdot\eta-v\cdot\xi\big)\right) \V^{\alpha\beta}(u-\xi,v-\eta,t-\tau) d\xi d\eta d\tau
\\\\ \ds~~~~~~~~~~~~~~~~~~
~=~\iiint_{\R^{2n+1}} f\left(\xi,\eta,\tau+\mu\big(u\cdot\eta-v\cdot\xi\big)\right) 
\\ \ds~~~~~~~~~~~~~~~~~~~~~~~~~~
|u-\xi|^{\alpha-n}|v-\eta|^{\alpha-n}|t-\tau|^{\beta-1} \Bigg[ {|u-\xi||v-\eta|\over |t-\tau|}+{|t-\tau|\over |u-\xi||v-\eta|}\Bigg]^{-{|\alpha-n\beta |\over n+1}}d\xi d\eta d\tau.
\end{array}
\eeq
Consider a more general situation by replacing $\V^{\alpha\beta}(\xi,\eta,\tau)$  with
\bel{V gneral}
|\xi|^{\alpha_1-n}|\eta|^{\alpha_2-n}|\tau|^{\beta-1} \Bigg[ {|\xi||\eta|\over |\tau|}+{|\tau|\over |\xi||\eta|}\Bigg]^{-\vartheta},\qquad \alpha_1,\alpha_2,\beta\in\R,\qquad \vartheta>0.
\eeq
By changing dilations $(u,v,t)\mt (\delta_1 u, \delta_2 v, \delta_1\delta_2\lambda t)$ and $(\xi,\eta,\tau)\mt (\delta_1 \xi, \delta_2 \eta, \delta_1\delta_2\lambda \tau)$ for $\delta_1,\delta_2>0$ and  $0<\lambda<1$ or $\lambda>1$,  we have
\bel{Dila I general}
\begin{array}{lr}\ds
\left\{ \iiint_{\R^{2n+1}} 
\left\{\iiint_{\R^{2n+1}} f\left[\delta^{-1}_1 \xi,\delta^{-1}_2\eta,\delta^{-1}_1\delta^{-1}_2\lambda^{-1}\big[\tau+\mu\lambda\big(u\cdot\eta-v\cdot\xi\big)\big]\right] \right.\right.
\\ \ds
\left. \left.|u-\xi|^{\alpha_1-n}|v-\eta|^{\alpha_2-n}|t-\tau|^{\beta-1} \Bigg[ {|u-\xi||v-\eta|\over |t-\tau|}+{|t-\tau|\over |u-\xi||v-\eta|}\Bigg]^{-\vartheta}
 d\xi d\eta d\tau\right\}^q dudvdt\right\}^{1\over q}
\\\\ \ds
=~\delta_1^{\alpha_1+\beta}\delta_2^{\alpha_2+\beta}\delta_1^{n+1\over q}\delta_2^{n+1\over q}\lambda^\beta\lambda^{1\over q}~\left\{ \iiint_{\R^{2n+1}} 
\left\{\iiint_{\R^{2n+1}} f\left( \xi,\eta,\tau+\mu\big(u\cdot\eta-v\cdot\xi\big)\right) \right.\right.
\\ \ds~~~
\left. \left.|u-\xi|^{\alpha_1-n}|v-\eta|^{\alpha_2-n}|t-\tau|^{\beta-1} \Bigg[ {|u-\xi||v-\eta|\over \lambda|t-\tau|}+{\lambda|t-\tau|\over |u-\xi||v-\eta|}\Bigg]^{-\vartheta}
 d\xi d\eta d\tau\right\}^q dudvdt\right\}^{1\over q}
\\\\ \ds
\ge~\delta_1^{\alpha_1+\beta}\delta_2^{\alpha_2+\beta}\delta_1^{n+1\over q}\delta_2^{n+1\over q}\lambda^\beta\lambda^{1\over q}\left\{ \begin{array}{lr}\ds \lambda^\vartheta,\qquad0<\lambda<1,
\\ \ds
\lambda^{-\vartheta},\qquad \lambda>1
\end{array}\right.
\\ \ds~~~
\left\{ \iiint_{\R^{2n+1}} 
\left\{\iiint_{\R^{2n+1}} f\left( \xi,\eta,\tau+\mu\big(u\cdot\eta-v\cdot\xi\big)\right) \right.\right.
\\ \ds~~~
\left. \left.|u-\xi|^{\alpha_1-n}|v-\eta|^{\alpha_2-n}|t-\tau|^{\beta-1} \Bigg[ {|u-\xi||v-\eta|\over |t-\tau|}+{|t-\tau|\over |u-\xi||v-\eta|}\Bigg]^{-\vartheta}
 d\xi d\eta d\tau\right\}^q dudvdt\right\}^{1\over q}.
\end{array}
\eeq
The $\L^p\mt\L^q$-norm inequality in (\ref{Result Two}) implies that the last line of (\ref{Dila I general}) is bounded by
\bel{I L^p}
\begin{array}{lr}\ds
\left\{ \iiint_{\R^{2n+1}} \Big[ f\left(\delta^{-1}_1 \xi,\delta^{-1}_2\eta,\delta^{-1}_1\delta^{-1}_2\lambda^{-1} \tau\right) \Big]^p d\xi d\eta d\tau\right\}^{1\over p}~=~\delta_1^{n+1\over p}\delta_2^{n+1\over p} \lambda^{1\over p} \left\| f\right\|_{\L^p(\R^{2n+1})}.
\end{array}
\eeq
Again, this must be true for every $\delta_1,\delta_2>0$ and $0<\lambda<1$ or $\lambda>1$. We necessarily have
\bel{Constraints alphabeta}
\begin{array}{cc}\ds
{\alpha_1+\beta\over n+1}~=~{1\over p}-{1\over q}~=~{\alpha_2+\beta\over n+1},
\\\\ \ds
\beta+\vartheta~\ge~{1\over p}-{1\over q}\qquad\hbox{\small{or}}\qquad \beta-\vartheta~\leq~{1\over p}-{1\over q}.
\end{array}
\eeq
The first constraint in (\ref{Constraints alphabeta}) forces us to have $\alpha_1=\alpha_2$. Therefore, write
\bel{Formula 12}
{\alpha+\beta\over n+1}~=~{1\over p}-{1\over q}
\eeq
where $\alpha=\alpha_1=\alpha_2$.
By bringing this to the second constraint in (\ref{Constraints alphabeta}), we find
\bel{theta computa}
\begin{array}{lr}\ds
\vartheta~\ge~\beta-{\alpha+\beta\over n+1}~=~{n\beta-\alpha\over n+1}\qquad \hbox{\small{or}}
\qquad
\vartheta~\ge~{\alpha+\beta\over n+1}-\beta~=~{\alpha-n\beta\over n+1}.
\end{array}
\eeq
Together, we conclude
\bel{theta constraint}
\vartheta~\ge~{|\alpha-n\beta|\over n+1}.
\eeq

\section{Size comparison between kernels}
\setcounter{equation}{0}
Let $0<\a<n$, $0<\b<1$ and $\a\ge n\b $.  We aim to show
\bel{size compara}
\begin{array}{lr}\ds
\Omega^{\a\b}(\xi,\eta,\tau)~=~\left[\frac{1}{|\xi|^2+|\eta|^2}\right]^{n-\a}\left[\frac{1}{|\xi|^2+|\eta|^2+|\tau|}\right]^{1-\b}
\\\\ \ds~~~~~~~~~~~~~~~~~~~
~\leq~|\xi|^{\alpha-n}|\eta|^{\alpha-n}|\tau|^{\beta-1}\left[\frac{|\xi||\eta|}{|\tau|}+\frac{|\tau|}{|\xi||\eta|}\right]^{-\frac{|\alpha-n\beta |}{n+1}}
\end{array}
\eeq
for some $\alpha,\beta$ satisfying $\alpha+\beta=\a+\b$.

Observe that
\bel{size}
\Omega^{\a\b}(\xi,\eta,\tau)
~\leq~\left(\frac{1}{|\xi||\eta|}\right)^{n-\a}\left[\frac{1}{|\xi||\eta|+|\tau|}\right]^{1-\b}.
\eeq
Suppose $|\xi||\eta|\geq|\tau|$. We further bound (\ref{size}) by
\bel{uv>t}
\begin{array}{lr}\ds
|\xi|^{\a-n}|\eta|^{\a-n}|\tau|^{\b-1}\left[\frac{|\tau|}{|\xi||\eta|+|\tau|}\right]^{1-\b}
~\leq~|\xi|^{\a-n}|\eta|^{\a-n}|\tau|^{\b-1}\left[\frac{|\xi||\eta|}{|\tau|}+\frac{|\tau|}{|\xi|\eta|}\right]^{-(1-\b)}.
\end{array}
\eeq
Choose 
\bel{ab case1}
\alpha~=~\a,\qquad \beta~=~\b.
\eeq
 We find
\bel{conclusion 1}
1-\b-\frac{|\alpha-n\beta|}{n+1}~=~1-\b-\frac{\a-n\b}{n+1}~=~1-\frac{\a+\b}{n+1}~>~0.
\eeq
Combining (\ref{size})-(\ref{uv>t}) and (\ref{ab case1})-(\ref{conclusion 1}) gives us 
\bel{boundedness uv>t}
\begin{array}{lr}\ds
\Omega^{\a\b}(\xi,\eta,\tau)~\leq~|\xi|^{\a-n}|\eta|^{\a-n}|\tau|^{\b-1}\left[\frac{|\xi||\eta|}{|\tau|}+\frac{|\tau|}{|\xi|\eta|}\right]^{-(1-\b)}
\\\\ \ds~~~~~~~~~~~~~~~~~~~
~\leq~|\xi|^{\alpha-n}|\eta|^{\alpha-n}|\tau|^{\beta-1}\left[\frac{|\xi||\eta|}{|\tau|}+\frac{|\tau|}{|\xi||\eta|}\right]^{-\frac{|\alpha-n\beta|}{n+1}}.
\end{array}
\eeq
Suppose $|\xi||\eta|\leq|\tau|$. Assert  $\b<\theta<1-\frac{\a-n\b}{n+1}$. We further bound (\ref{size})  by
\bel{uv<t}
\begin{array}{lr}\ds
|\xi|^{\a-n}|\eta|^{\a-n}\left[\frac{1}{|\xi||\eta|+|\tau|}\right]^{\theta-\b}\left[\frac{1}{|\xi||\eta|+|\tau|}\right]^{1-\theta}
\\\\ \ds
\leq~|\xi|^{\a-n}|\eta|^{\a-n}\left(\frac{1}{|\tau|}\right)^{\theta-\b}\left[\frac{1}{|\xi||\eta|+|\tau|}\right]^{1-\theta}
~=~|\xi|^{\a-n+\theta-1}|\eta|^{\a-n+\theta-1}|\tau|^{\b-\theta}\left[\frac{|\xi||\eta|}{|\xi||\eta|+|\tau|}\right]^{1-\theta}
\\\\ \ds
\leq~|\xi|^{(\a+\theta-1)-n}|\eta|^{(\a+\theta-1)-n}|\tau|^{(\b-\theta+1)-1}\left[\frac{|\xi||\eta|}{|\tau|}+\frac{|\tau|}{|\xi|\eta|}\right]^{-(1-\theta)}.
\end{array}
\eeq
Choose 
\bel{ab case2}
\alpha~=~\a+\theta-1,\qquad \beta~=~\b-\theta+1.
\eeq 
Because $\theta<1-\frac{\a-n\b}{n+1}$, we have
\bel{conclusion 2}
\begin{array}{lr}\ds
1-\theta-\frac{\left|\alpha-n\beta\right|}{n+1}~=~1-\theta-\frac{\left|\a-n\b+(n+1)\theta-(n+1)\right|}{n+1}
\\\\ \ds~~~~~~~~~~~~~~~~~~~~~~~~~~~
~=~1-\theta+\frac{\a-n\b}{n+1}+\theta-1~=~\frac{\a-n\b}{n+1}~\ge~0.
\end{array}
\eeq
By putting together (\ref{size}), (\ref{uv<t}) and (\ref{ab case2})-(\ref{conclusion 2}),  we obtain
\bel{boundedness uv<t}
\begin{array}{lr}\ds
\Omega^{\a\b}(\xi,\eta,\tau)~\leq~|\xi|^{(\a+\theta-1)-n}|\eta|^{(\a+\theta-1)-n}|\tau|^{(\b-\theta+1)-1}\left[\frac{|\xi||\eta|}{|\tau|}+\frac{|\tau|}{|\xi|\eta|}\right]^{-(1-\theta)}
\\\\ \ds~~~~~~~~~~~~~~~~~~~
~\leq~|\xi|^{\alpha-n}|\eta|^{\alpha-n}|\tau|^{\beta-1}\left[\frac{|\xi||\eta|}{|\tau|}+\frac{|\tau|}{|\xi||\eta|}\right]^{-\frac{|\alpha-n\beta|}{n+1}}.
\end{array}
\eeq

\section{Proof of Theorem Two}
\setcounter{equation}{0}
Let
\bel{Formula}
{\alpha+\beta\over n+1}~=~{1\over p}-{1\over q},\qquad 1<p<q<\infty
\eeq
which is an necessity  for the $\L^p\mt\L^q$-norm inequality in (\ref{Result Two}).  

We now turn to prove the converse. First, as shown in (\ref{V}),  $\V^{\alpha\beta}$  is positive definite. Therefore,  it is suffice to assert $f\ge0$.

Suppose  $\alpha\ge n\beta $. We have ${|\alpha-n\beta |\over n+1}={\alpha-n\beta \over n+1}$ and
\bel{EST1}
\begin{array}{lr}\ds
\V^{\alpha\beta}(\xi,\eta,\tau)~=~|\xi|^{\alpha-n}|\eta|^{\alpha-n}|\tau|^{\beta-1} \Bigg[ {|\xi||\eta|\over |\tau|}+{|\tau|\over |\xi||\eta|}\Bigg]^{-{\alpha-n\beta \over n+1}}
\\\\ \ds~~~~~~~~~~~~~~~~~~~
~\leq~|\xi|^{\alpha-n}|\eta|^{\alpha-n}|\tau|^{\beta-1} \Bigg[ {|\xi||\eta|\over |\tau|}\Bigg]^{-{\alpha-n\beta \over n+1}}
~=~|\xi|^{n\big[{\alpha+\beta\over n+1}\big]-n}|\eta|^{n\big[{\alpha+\beta\over n+1}\big]-n} |\tau|^{{\alpha+\beta\over n+1}-1}
\end{array}
\eeq
for $\xi\neq0$, $\eta\neq 0$,   $\tau\neq0$.

Suppose $\alpha\leq n\beta $. We find  ${|\alpha-n\beta |\over n+1}={n\beta -\alpha\over n+1}$ and
\bel{EST2}
\begin{array}{lr}\ds
\V^{\alpha\beta}(\xi,\eta,\tau)~=~|\xi|^{\alpha-n}|\eta|^{\alpha-n}|\tau|^{\beta-1} \Bigg[ {|\xi||\eta|\over |\tau|}+{|\tau|\over |\xi||\eta|}\Bigg]^{{\alpha-n\beta \over n+1}}
\\\\ \ds~~~~~~~~~~~~~~~~~~
~\leq~|\xi|^{\alpha-n}|\eta|^{\alpha-n}|\tau|^{\beta-1} \Bigg[ {|\tau|\over |\xi||\eta|}\Bigg]^{\alpha-n\beta \over n+1}
~=~|\xi|^{n\big[{\alpha+\beta\over n+1}\big]-n}|\eta|^{n\big[{\alpha+\beta\over n+1}\big]-n} |\tau|^{{\alpha+\beta\over n+1}-1} 
\end{array}
\eeq
 for $\xi\neq0$, $\eta\neq 0$,   $\tau\neq0$.

 Let $\I_{\alpha\beta}$ defined in (\ref{V})-(\ref{I}).
By changing variable $\tau\mt \tau+\mu\big(u\cdot\eta-v\cdot\xi\big)$, we have 
 \bel{I < =}
\begin{array}{lr}\ds
\I_{\alpha\beta} f(u,v,t)~=~\iiint_{\R^{2n+1}} f\left(\xi,\eta,\tau+\mu\big(u\cdot\eta-v\cdot\xi\big)\right) 
\\ \ds~~~~~~~~~~~~~~~~~~~~~~~~~~
|u-\xi|^{\alpha-n}|v-\eta|^{\alpha-n}|t-\tau|^{\beta-1} \Bigg[ {|u-\xi||v-\eta|\over |t-\tau|}+{|t-\tau|\over |u-\xi||v-\eta|}\Bigg]^{-{|\alpha-n\beta |\over n+1}}d\xi d\eta d\tau
\\\\ \ds~~~~~~~~~~~~~~~~~~
~\leq~
\iiint_{\R^{2n+1}}  f\left(\xi,\eta,\tau+\mu\big(u\cdot\eta-v\cdot\xi\big)\right) 
\\ \ds~~~~~~~~~~~~~~~~~~~~~~~~~~
|u-\xi|^{n\big[{\alpha+\beta\over n+1}\big]-n}|v-\eta|^{n\big[{\alpha+\beta\over n+1}\big]-n} |t-\tau|^{{\alpha+\beta\over n+1}-1}d\xi d\eta d\tau\qquad\hbox{\small{by (\ref{EST1})-(\ref{EST2})}}.
\end{array}
\eeq
 Define
\bel{F_beta}
\F_{\alpha\beta}(\xi,\eta, u,v,t)~=~\int_\R f\left(\xi,\eta,\tau+\mu\left(u \cdot\eta-v \cdot\xi\right)\right) |t-\tau|^{{\alpha+\beta\over n+1}-1} d\tau.
\eeq
From (\ref{I < =})-(\ref{F_beta}), we find
\bel{I rewrite F}
\begin{array}{lr}\ds
\I_{\alpha\beta} f(u,v,t)
~\leq~\iint_{\R^{2n}} |u-\xi|^{n\big[{\alpha+\beta\over n+1}\big]-n}|v-\eta|^{n\big[{\alpha+\beta\over n+1}\big]-n} \F_{\alpha\beta}(\xi,\eta, u,v,t)d\xi d\eta.
\end{array}
\eeq
Recall  the {\bf Hardy-Littlewood-Sobolev theorem} stated in the beginning of this paper. By applying (\ref{HLS Ineq}) with $\a={\alpha+\beta\over n+1}$ and $\N=1$, we have
\bel{F regularity}
\begin{array}{lr}\ds
\left\{\int_{\R} \F^q_{\alpha\beta}(\xi,\eta, u,v,t) dt \right\}^{1\over q}~\leq~\C_{p~q} \left\{\int_{\R} \Big[ f\left(\xi,\eta,t+\mu\left(u \cdot\eta-v \cdot\xi\right)\right)\Big]^p dt\right\}^{1\over p}
\\\\ \ds~~~~~~~~~~~~~~~~~~~~~~~~~~~~~~~~~~~~~~~~~
~=~\C_{p~q} \left\| f(\xi,\eta,\cdot)\right\|_{\L^p(\R)}
\end{array}
\eeq
regardless of $(u,v)\in\R^n\times\R^n$.

On the other hand, by applying (\ref{HLS Ineq}) with $\a=n \big[{\alpha+\beta\over n+1}\big]$ and $\N=n$, we find
\bel{n Regularity}
\begin{array}{lr}\ds
\left\{\int_{\R^n}         \left\{ \int_{\R^n}   |u-\xi|^{n \big[{\alpha+\beta\over n+1}\big]-n}  \left\| f(\xi,\eta,\cdot)\right\|_{\L^p(\R)} d\xi\right\}^q du\right\}^{1\over q}~\leq~\C_{p~q} ~ \left\{\int_{\R^n} \left\| f(u,\eta,\cdot)\right\|_{\L^p(\R)}^p du\right\}^{1\over p},
\\\\ \ds
\left\{\int_{\R^n}         \left\{ \int_{\R^n}   |v-\eta|^{n \big[{\alpha+\beta\over n+1}\big]-n}  \left\| f(\xi,\eta,\cdot)\right\|_{\L^p(\R)} d\eta\right\}^q dv\right\}^{1\over q}~\leq~\C_{p~q} ~ \left\{\int_{\R^n} \left\| f(\xi,v,\cdot)\right\|_{\L^p(\R)}^p dv\right\}^{1\over p}.
\end{array}
\eeq

From (\ref{I rewrite F}), we have
\bel{EST}
\begin{array}{lr}\ds
\left\|\I_{\alpha\beta} f\right\|_{\L^q(\R^{2n+1})}
\\\\ \ds
~\leq~\left\{ \iiint_{\R^{2n+1}}\left\{\iint_{\R^{2n}} |u-\xi|^{n \big[{\alpha+\beta\over n+1}\big]-n}|v-\eta|^{n \big[{\alpha+\beta\over n+1}\big]-n} \F_{\alpha\beta}(\xi,\eta, u,v,t)d\xi d\eta\right\}^q dudvdt\right\}^{1\over q}
\\\\ \ds
~\leq~      \left\{\iint_{\R^{2n}}      \left\{ \iint_{\R^{2n}}   |u-\xi|^{n \big[{\alpha+\beta\over n+1}\big]-n}|v-\eta|^{n \big[{\alpha+\beta\over n+1}\big]-n}  \left\{\int_\R \F^q_{\alpha\beta}(\xi,\eta, u,v,t) dt \right\}^{1\over q} d\xi d\eta\right\}^q dudv\right\}^{1\over q}
\\ \ds~~~~~~~~~~~~~~~~~~~~~~~~~~~~~~~~~~~~~~~~~~~~~~~~~~~~~~~~~~~~~~~~~~~~~~~~~~~~~~~~~~~~~~~~~
\hbox{\small{by Minkowski integral inequality}}
\\\\ \ds
~\leq~\C_{p~q} \left\{\iint_{\R^{2n}}      \left\{ \iint_{\R^{2n}}   |u-\xi|^{n \big[{\alpha+\beta\over n+1}\big]-n}|v-\eta|^{n \big[{\alpha+\beta\over n+1}\big]-n}   \left\| f(\xi,\eta,\cdot)\right\|_{\L^p(\R)} d\xi d\eta\right\}^q dudv\right\}^{1\over q}\qquad\hbox{\small{by (\ref{F regularity})}}
\\\\ \ds
~\leq~\C_{p~q} \left\{\int_{\R^n}        \left\{\int_{\R^n}  \left\{ \int_{\R^n}   |u-\xi|^{\alpha-n}  \left\| f(\xi,v,\cdot)\right\|_{\L^p(\R)} d\xi\right\}^p dv\right\}^{q\over p} du\right\}^{1\over q}
\qquad
\hbox{\small{by  (\ref{n Regularity})}}
\\\\ \ds
~\leq~\C_{p~q} \left\{\int_{\R^n} \left\{\int_{\R^n}         \left\{ \int_{\R^n}   |u-\xi|^{\alpha-n}  \left\| f(\xi,v,\cdot)\right\|_{\L^p(\R)} d\xi\right\}^q du \right\}^{p\over q} dv\right\}^{1\over p}
\\ \ds~~~~~~~~~~~~~~~~~~~~~~~~~~~~~~~~~~~~~~~~~~~~~~~~~~~~
\hbox{\small{by Minkowski integral inequality}}
\\\\ \ds
~\leq~\C_{p~q}      \left\{\iint_{\R^{2n}}  \left\| f(u,v,\cdot)\right\|_{\L^p(\R)}^p du dv\right\}^{1\over p}\qquad\hbox{\small{by  (\ref{n Regularity})}}
\\\\ \ds
~=~\C_{p~q} \left\| f\right\|_{\L^p(\R^{2n+1})}.
\end{array}
\eeq

{\small Westlake University}\\
{\small sunchuhan@westlake.edu.cn}

{\small  Westlake University}\\
{\small wangzipeng@westlake.edu.cn}

\end{document}